\documentclass[journal,onecolumn,12pt]{IEEEtran}
\usepackage{amsmath} 
\usepackage{amsfonts}
\usepackage{amssymb}
\usepackage{amsthm}
\usepackage{makecell}
\usepackage{bbm}
\usepackage[margin=0.75in]{geometry} 
\usepackage{xcolor}
\newtheorem{theorem}{Theorem}

\newtheorem{proposition}{Proposition}

\renewcommand{\l}{\ell}
\newcommand{\Tr}{\mathrm{Tr}}
\newcommand{\E}{\mathbb{E}}
\newcommand{\F}{\mathbb{F}}
\newcommand{\C}{\mathcal{C}}

\title{Comments on: ``Symmetric Pseudo-Random Matrices''}
\author{Chin Hei Chan\thanks{C. H. Chan is with Hetao Institute of Mathematics and Interdisciplinary Sciences (HIMIS), Shenzhen 518017, Guangdong, China (email address: chenzhanxi@himis-sz.cn). He is supported by the research start-up fund from HIMIS.}}
\date{}
\begin{document}
\maketitle
\begin{abstract}
    In 2018 \cite{SPRM}, Soloveychik, Xiang and Tarokh considered a pseudo-random symmetric circulant matrix constructed from binary Golomb sequences of length $n=2^m-1$, and claimed a proof that its empirical spectral distribution converges almost surely to the semicircle law as $n$ grows to infinity by the method of moments. In this comment note we show that their argument contains several technical flaws and inappropriate applications of technical lemmas. Instead we demonstrate that the eigenvalues of the matrix are simply a normalized twisted Kloosterman sum varying over the multiplicative character, to which we apply Katz's result \cite{Katz} directly to establish the asymptotically semicircle spectral distribution deterministically.
\end{abstract}
\section{Introduction}
In 2018 \cite{SPRM}, Soloveychik, Xiang and Tarokh considered a random matrix constructed as follows:

Let $\C$ be the $[n=2^m-1, m]_2$ simplex code constructed from a primitive binary polynomial $f$ of degree $m$. Pick a nonzero codeword $\varphi \in \C$ (such a codeword is called a Golomb sequence or m-sequence). Let $\mathbf{T}$ be the binary circulant $n \times n$ matrix defined by
$$\mathbf{T}=(t_{ij})_{i,j=0}^{n-1}=(\varphi(j-i))_{i,j=0}^{n-1}$$
where the input of $\varphi$ is interpreted modulo $n$.

In particular, the first row of $\mathbf{T}$ is simply the codeword $\varphi$ and the rest rows are successive cyclic shifts of $\varphi$, which are in fact the rest nonzero codewords of $\mathcal{C}$.

Now define the real $n \times n$ matrix
$$\mathbf{A}_n=(a_{ij})_{i,j=0}^{n-1}=\frac{1}{2\sqrt{n}}((-1)^{t_{ij}+t_{ji}})_{i,j=0}^{n-1}=\frac{1}{2\sqrt{n}}\zeta(\mathbf{T}+\mathbf{T}^T),$$
where $\zeta(x)=(-1)^x$ for $x \in \F_2$.

They consider the ensemble $\mathcal{A}_n$ of all matrices
$$\mathbf{A}_n(a):=\frac{1}{2\sqrt{n}}((-1)^{\varphi(j-i+a)+\varphi(i-j+a)})_{i,j=0}^{n-1}$$
for $0 \leq a \leq n-1$ and their negatives, endowed with the uniform probability measure.

They performed simulation experiments on 100 instants of the ensemble with $m=13,14$ respectively, which suggested that the average empirical spectral distribution of $\mathbf{A}_n(a)$ converges to the semicircle law with probability density
$$f_\mathrm{SC}(x)=\frac{2}{\pi}\sqrt{1-x^2}, -1 \leq x \leq 1.$$

Then in the \textbf{Appendix} section, they derived a mathematical proof using the method of moments. In particular, they claimed the following two results, where $\beta_r(\mathbf{A}_n)$ denotes the $r$-th spectral moment of $\mathbf{A}_n$, $\E[\cdot]$ and $\mathbb{V}[\cdot]$ denote expectation and variance over $\mathcal{A}_n$ respectively:
\begin{proposition}\label{P1}\cite[Proposition 1]{SPRM}
    Let $\mathbf{A}_n \in \mathcal{A}_n$. For a fixed $r \in \mathbb{N}$ and $n=n(m) \to \infty$,
    $$\E[\beta_r(\mathbf{A}_n)]=\begin{cases}
        \frac{C_{r/2}}{2^r}+O\left(\frac{1}{n}\right), &r \text{ even}\\
        0, &r \text{ odd}.
    \end{cases}$$
\end{proposition}
\begin{proposition}\label{P2}\cite[Proposition 2]{SPRM}
    Let $\mathbf{A}_n \in \mathcal{A}_n$. For a fixed $r \in \mathbb{N}$ and $n=n(m) \to \infty$,
    $$\mathbb{V}[\beta_r(\mathbf{A}_n)]=O\left(\frac{1}{n^2}\right).$$
\end{proposition}
Upon thorough examination of their proof, we find that there are several technical flaws and inappropriate applications of technical lemmas without sufficient explanation. In particular, we demonstrate that Proposition \ref{P2} is false at least for $r=1$.

On the other hand, we discover that the eigenvalues of the matrix can be expressed in terms of a normalized twisted Kloosterman sum, whose equidistribution for the semicircle measure was already mathematically verified by Katz in 2012 \cite{Katz}. This in fact establishes the stronger conclusion stated in \cite[Remark 1]{SPRM} that the empirical spectra of the sequence of matrices $\{\mathbf{A}_n\}$ converge to the semicircle law when $n$ grows, in a deterministic sense.

The comment note is organized as follows. In Section \ref{TF} we provide a list of technical flaws and inappropriate applications of technical lemmas discovered, with concrete explanations. In Section \ref{KS} we provide the derivation of the eigenvalues as a normalized twisted Kloosterman sum that has been studied in \cite{Katz}.

\section{Technical Flaws}\label{TF}
In this section, we state all the technical flaws and inappropriate applications of technical lemmas in \cite{SPRM}.
\subsection{Computation Error}
\begin{enumerate}
    \item \emph{The derivation of Eq. (76) from Eq. (67)}
    
        Eq. (76) is derived from Eq. (67) by converting $i_q$'s to $t_q$'s, where $t_q$ is defined as in Eq. (69). Comparing the two equations, this amounts to claiming
    \begin{equation}\label{E1}\tag{*}
    \sum_{i_0,\cdots,i_{r-1}=0}^{n-1}(-1)^{\sum_{q=0}^{r-1}(\varphi(i_{q+1}-i_q+a)+\varphi(i_q-i_{q+1}+a))}=\sum_{\mathbf{t} \in \mathcal{T}_r}(-1)^{\tau(\nu(\mathbf{t});a)},
    \end{equation}
    where $\mathcal{T}_r,\nu(\mathbf{t})$ and $\tau(\nu(\mathbf{t});a)$ are defined in Eqs. (73), (70) and (74) respectively.

    While the congruence $\sum_{q=0}^{r-1}(\varphi(i_{q+1}-i_q+a)+\varphi(i_q-i_{q+1}+a)) \equiv \tau(\nu(\mathbf{t});a) \pmod{2}$ follows directly from Eqs. (68) and (74), the authors did not realize that the summation sets $\{(i_0,i_1,\cdots,i_{r-1}): i_q \in [n]\}$ and $\mathcal{T}_r$ are not of the same cardinality (the former is $n^r$ but the latter is $n^{r-1}$). In fact the mapping (Eq. (68)) between these two sets is $n$-to-one. Indeed, adding the same integer modulo $n$ to all $i_q$'s preserves $\mathbf{t}$. Hence there is a missing factor of $n$ on the RHS of (\ref{E1}).
    
\end{enumerate}
\subsection{False Claims}
\begin{enumerate}
    \item \emph{Non-even cycles versus $\nu(\mathbf{t}) \neq \mathbf{0}$}
    
    In computing III (just after Eq. (78)), they simply assumed that non-even cycles yield $\nu(\mathbf{t}) \neq \mathbf{0}$ without justification. In fact a counterexample can be easily constructed: let $r=4$ and $i_0=1,i_1=2,i_2=4,i_3=3$. All the edges in the cycle ($(1,2),(2,4),(4,3),(3,1)$) have multiplicity one (hence the cycle is non-even). On the other hand, by Eq. (68), $t_0=1,t_1=2,t_2=-1,t_3=-2$. Hence $\nu(\mathbf{t})=\mathbf{0}$ by Eq. (70).

    This implies that the set of $\mathbf{t}$'s with $\nu(\mathbf{t})=\mathbf{0}$ should be larger than the set of even cycles. Indeed, if $t_q$'s come in $\pm i \pmod{n}$ pairs, then both conditions $\nu(\mathbf{t})=\mathbf{0}$ and $\sum_{q=0}^{n-1} t_q=0$ satisfy trivially. This gives $(r-1)!!n^{r/2}+O(n^{r/2-1})$ possibilities for $\mathbf{t}$, and it can be shown that this class contributes precisely to the highest-order term. Note that the numbers $(r-1)!!$ are exactly the $r$-th moments of the standard normal distribution when $r$ is even. This is consistent with the fact that the expected empirical spectral distribution of a truly random symmetric circulant matrix converges to the normal distribution, in which the authors also knew (see \cite[Section VII-B]{SPRM}).

    With this modification where the true II (corresponds to $\nu(\mathbf{t})=\mathbf{0}$) gives the Gaussian moments $\frac{(r-1)!!}{2^r}$ as the leading term when $r$ is even, if Proposition \ref{P1} does hold, then the true III (corresponding to $\nu(\mathbf{t}) \neq \mathbf{0}$) should converge to $\frac{C_{r/2}-(r-1)!!}{2^r}$, which is negative but not negligible for $r \geq 4$. The incorrectly claimed order ($O(n^{-1})$ instead of $O(1)$) can be partially explained by the missing $n$ factor of (\ref{E1}), but it turns out that there are other flaws in the computation of III, which will be discussed in the forthcoming paragraphs.
   
     \item \emph{Upper bound for $\rho_r(\l)$}
     
     In Eq. (86), they claimed a uniform upper bound for $\rho_r(\l)$ independent of $n$. This is not valid when $\l < r$. Actually one can think in this way: for each $q$ in the first half of the support of the specified codeword, we assign either $q$ or $-q$ to some particular coordinate of $\mathbf{t}$. Then we have $\nu(\mathbf{t}')=\mathbf{0}$ where $\mathbf{t}'$ is $\mathbf{t}$ restricted to the unassigned coordinates. Now it is easy to see that there are $n$-polynomially many choices for $\mathbf{t}'$. As a further remark, under the restricted condition (72), the quantity $\rho'_r(\mathbf{c})$ that counts the number of $\mathbf{t} \in \mathcal{T}_r$ producing the same codeword $\mathbf{c}$ actually depends on $\mathbf{c}$ itself but \emph{not only} its weight anymore, so the idea of deriving III from III' in the paper is also invalid.

      \item \emph{Weights of $\mathcal{H}$}
    
    After Eq. (80), they claimed that the weight of any codeword in $\mathcal{H}$ must be divisible by 4. This is not true in general. They overlooked the fact that the zeroth bit of $\nu(\mathbf{t})$ is always 0 (or even) regardless of the number of zeros in $\mathbf{t}$, since zeros are double-counted in the definition of $\nu(\mathbf{t})$. In fact one can only conclude that the weight of any nonzero codeword in $\mathcal{H}$ is even and at least 4 (the latter is due to the fact that the linear code $\mathcal{C}^\bot$, which is a Hamming code, has minimum distance 3). Nevertheless, this is a minor flaw that does not influence the forthcoming argument significantly.
    \end{enumerate}
    \subsection{Inappropriate Applications of Technical Lemmas}
    \begin{enumerate}
    \item \emph{Applying Lemma 7 to $\langle\overline{\mathcal{H}'}\rangle$}
    
    In Eq. (98), they applied Lemma 7 to $\langle\overline{\mathcal{H}'}\rangle$ directly without any justification. Lemma 7 applies to codes under the conditions of Lemma 6, that is, binary linear codes of the form $S(D,G)$ where the divisors $D,G$ satisfy certain conditions. We verified that $\langle \overline{\mathcal{H}'}\rangle$ is only \emph{related to the half-puncturing} of some $S(D,G)$. Moreover, it happens that the corresponding $G_-$ is a divisor of degree 2 on a single point rather than just a single point as required in Lemma 6, so Lemma 7 cannot be applied directly here. 
    
    \item \emph{Applying Lemma 6 to a misinterpreted multiplicative condition}
    
    After Eq. (100), they claimed that $\alpha$ in Lemma 6 increases by 1 under the linear relation Eq. (72). However this relation is \emph{not} linear over $\F_2$ or $\F_{2^m}$ (but rather multiplicative, under the isomorphism $\mathbb{Z}/n\mathbb{Z} \to \F_{2^m}^\times, t \mapsto \beta^t$ for a fixed primitive $\beta \in \F_{2^m}^\times$), so Lemma 6 is not applicable. This presents another obstacle to derive III from III'.
    \end{enumerate}
    \subsection{Order of the Variance in Proposition \ref{P2}}
      
      In Eq. (102), they made a mistake similar to that in Eq. (76). In fact since it involves a sum over two sets of indices, the RHS of Eq. (102) is smaller than it should be by a factor of $n^2$. Disregarding any further mistakes (for instance, they assumed that $U_0$ corresponds exactly to even paths), what they obtained for the variance is actually $O(1)$ only. Consequently their claim of almost sure convergence in \cite[Corollary 2]{SPRM} is not verified (nor is convergence in probability).
      
      In fact the variance of the first moment $\beta_1(\mathbf{A}_n)$ can be computed directly without any combinatorial argument. Indeed $\beta_1(\mathbf{A}_n)$ is simply any diagonal entry of the circulant matrix $\mathbf{A}_n$, which is a Rademacher variable divided by $2\sqrt{n}$; hence its variance is $\frac{1}{4n}$. This in particular disproves Proposition \ref{P2} for $r=1$.

\subsection{Definition of Standard Semicircle Law}
While in some existing literature the standard semicircle law is indeed defined as in \cite[Eq. (21)]{SPRM}, in parallel with the standard normal distribution, it is more convenient to define the \emph{standard} (or \emph{standardized}) semicircle law as the one with unit variance (this convention is also more common in random matrix theory). In this case, the probability density function is given by
\begin{equation}\label{SC}
    \varrho_\mathrm{SC}(x)=\frac{1}{2\pi}\sqrt{4-x^2}
\end{equation}
over $[-2,2]$, whose associated random variable is twice that defined in \cite[Eq. (21)]{SPRM}.

\section{Amended Proof}\label{KS}
\subsection{Problem Set-up}
In our setting, we use an algebraic representation of the Golomb sequences. Indeed, the nonzero codewords of $\mathcal{C}$ can be expressed as
$$\varphi_a=(\varphi_a(r))_{r=0}^{n-1}=(\Tr_1^m(a\alpha^r))_{r=0}^{n-1}$$
for $a \in \F_q^\times$, where $q=2^m=n+1$ and $\alpha$ is a fixed primitive element of $\F_q^\times$.

We consider the ensemble $\mathcal{A}_n'$ of matrices 
$$\mathbf{A}_n'=\mathbf{A}_n'(a)=(a_{ij}'(a))_{i,j=0}^{n-1}$$
and their negatives, where $a_{ij}'(a):=\frac{1}{\sqrt{n}}(-1)^{\varphi_a(j-i)+\varphi_a(i-j)}$, with the input of $\varphi_a$ interpreted modulo $n$. Here we use a different normalization $\frac{1}{\sqrt{n}}$ (compared with $\frac{1}{2\sqrt{n}}$ in \cite{SPRM}) to conform to the standard scale of probability distributions. 

We aim to prove the analogue of \cite[Remark 1]{SPRM}, namely:
\begin{theorem}\label{Thm}
    Let $\mathbf{A}_n' \in \mathcal{A}_n'$ for all $m$. Then the spectral measure of $\mathbf{A}_n'$ converges to the standard semicircle measure $\varrho_\mathrm{SC}$ as $m$ grows.
\end{theorem}
Here $\varrho_\mathrm{SC}$ is defined as in (\ref{SC}).

\subsection{Proof of Theorem \ref{Thm}}
Denote by $c_r(a):=\frac{1}{\sqrt{n}}(-1)^{\varphi_a(r)+\varphi_a(n-r)}$. Then it is easy to see that $c_{n-r}(a)=c_r(a)$ for all $r \in [1.. \frac{n-1}{2}]$, and $a_{ij}'(a)=c_{|i-j|}(a)$.

Instead of attempting to compute the expected spectral moments as in \cite{SPRM}, which requires a careful technical treatment of the III term, we exploit the symmetric circulant structure of $\pm A_n(a)$ to compute its eigenvalues directly as (see \cite[Remark 2]{Bose} for reference)
\begin{align}
\lambda_k^{\pm}(a)&=\sum_{r=0}^{n-1} \pm c_r(a)\cos\left(\frac{2\pi kr}{n}\right)\nonumber\\
&=\pm\sum_{r=0}^{n-1}\frac{(-1)^{\varphi_a(r)+\varphi_a(n-r)}}{\sqrt{n}}\cos\left(\frac{2\pi kr}{n}\right)\nonumber\\
&=\pm\frac{\sum_{r=0}^{n-1}\psi_a(\alpha^r+\alpha^{-r})\cos\left(\frac{2\pi kr}{n}\right)}{\sqrt{n}}\label{lambda},
\end{align}
where $\psi_a: \beta \mapsto (-1)^{\Tr_1^m(a\beta)}$ is a non-trivial additive character on $\F_q$.

Noting that $\cos\left(\frac{2\pi(n-k)r}{n}\right)=\cos\left(2r\pi-\frac{2\pi kr}{n}\right)=\cos\left(\frac{2\pi kr}{n}\right)$ for all $k$, we may rewrite the sum (\ref{lambda}) as
\begin{align}
    \lambda_k^\pm(a)&=\pm\frac{1+\sum_{r=1}^{(n-1)/2}\psi_a(\alpha^r+\alpha^{-r})\left(\cos\left(\frac{2\pi kr}{n}\right)+\cos\left(\frac{2\pi(n-k)r}{n}\right)\right)}{\sqrt{n}}\nonumber\\
    &=\pm\frac{1+\sum_{r=1}^{(n-1)/2}\psi_a(\alpha^r+\alpha^{-r})\left(2\cos\left(\frac{2\pi kr}{n}\right)\right)}{\sqrt{n}}\nonumber\\
    &=\pm\frac{1+\sum_{r=1}^{(n-1)/2}\psi_a(\alpha^r+\alpha^{-r})(2\Re(\chi_k(\alpha^r)))}{\sqrt{n}}\nonumber\\
    &=\pm\frac{\sum_{r=0}^{n-1}\psi_a(\alpha^r+\alpha^{-r})\chi_k(\alpha^r)}{\sqrt{n}}\label{lambda2},
\end{align}
where $\chi_k: \alpha \mapsto \exp\left(\frac{2\pi k\mathrm{i}}{n}\right)$ is a multiplicative character of $\F_q^\times$.

The numerator of (\ref{lambda2}) is precisely the twisted Kloosterman sum $K_{\psi_a}(1,1; \chi_k)$ (see \cite{Cohen}), where
$$K_\psi(t,u; \chi):=\sum_{\xi\in \F_q^\times} \psi(t\xi+u\xi^{-1})\chi(\xi)$$
for $t,u \in \F_q$ (not both zero), $\psi \in \widehat{(\F_q,+)} \setminus \{\psi_0\}$ and $\chi \in \widehat{\F_q^\times}$.

It happens that the sum $K_{\psi_a}(1,1; \chi_k)$ was first introduced by Ron Evans in May 2003 (strictly speaking his sum is $K_\psi(1,-1;\chi)$, but in characteristic 2 we have $1=-1$), who made empirical simulations of this sum for fixed $x$ and varying $\chi_k$ (equivalently, $k$), and conjectured that this sequence of $n$ sums, after normalization by $-\frac{1}{\sqrt{q}}$, is equidistributed for the standard semicircle measure (also known as the Sato-Tate measure) $\varrho_\mathrm{SC}$ in (\ref{SC}) on the closed interval $[-2,2]$, as $q \to \infty$. This was then proved by Katz in 2012 using techniques from algebraic geometry (see \cite[Chapter 14]{Katz} for the detailed proof).

Since
$$\lambda_k^\pm(a)=\mp\sqrt{1+\frac{1}{n}}\left(-\frac{1}{\sqrt{q}}K_{\psi_a}(1,1;\chi_k)\right)$$
and the semicircle density $\varrho_{\mathrm{SC}}(x)$ is symmetric over the $y$-axis, we conclude that both $\{\lambda_k^+(a)\}_{k=0}^{n-1}$ (eigenvalues of $\mathbf{A}_n'(a)$) and $\{\lambda_k^-(a)\}_{k=0}^{n-1}$ (eigenvalues of $-\mathbf{A}_n'(a)$) are equidistributed for $\varrho_{\mathrm{SC}}$ on $[-2,2]$ as $n \to \infty$ (equivalently $m \to \infty$).

This completes the proof of Theorem \ref{Thm}.

\section*{Acknowledgments}
The discovery of the technical flaws in \cite{SPRM} and that the eigenvalues of the matrix can be expressed as a twisted Kloosterman sum were initiated when the sole author was working in Department of Mathematics, Hong Kong University of Science and Technology. He would like to thank the department for providing financial support. Further literature review to discover Katz's proof on that sum and the drafting of the comment note are then completed in HIMIS.


\begin{thebibliography}{99}
\bibitem{Bose} A. Bose and J. Mitra, ``Limiting spectral distribution of a special circulant,'' \emph{Stat. Probab. Letters} \textbf{60} (2002), no. 1, 111--120.

\bibitem{Cohen} S. D. Cohen, ``Kloosterman sums and primitive elements in Galois fields,'' \emph{Acta Arithmetica} \textbf{94} (2000), no. 2, 173--201.


\bibitem{Katz} N. M. Katz, ``Convolution and Equidistribution: Sato-Tate Theorems for Finite-Field Mellin Transforms,'' \emph{Annals of Mathematics Studies} \textbf{180}. Princeton, NJ: Princeton University Press, 194 pp., 2012.


\bibitem{SPRM} I. Soloveychik, Y. Xiang and V. Tarokh, ``Symmetric Pseudo-Random Matrices,'' \emph{IEEE Trans. Inform. Theory} {\bf 64} (2018), no. 4, 3179--3196.


\end{thebibliography}
\end{document}